\documentclass[preprint,sort&compress]{elsarticle}
\usepackage{amssymb,amsthm,amsmath,caption,commath,bbm,psfrag,verbatim}
\usepackage{mathtools}
\usepackage{setspace,color,multirow,array,graphicx,epstopdf,float}
\usepackage{tikz}
\usetikzlibrary{arrows.meta, positioning}
\usepackage[FIGTOPCAP,nooneline]{subfigure}
\usepackage{capt-of}

\usepackage{tikz}
\usepackage{pgfplots}
\pgfplotsset{compat=1.18}
\usepackage{amsmath}

\usepgfplotslibrary{fillbetween}
\usepackage{tabularx}
\usepackage[margin=1in]{geometry}
\usepackage{tikz}
\usetikzlibrary{arrows.meta, positioning}
\usetikzlibrary{calc}
\usepackage{url} 
\usepackage{hyperref} 
\usepackage{newlfont}

\newtheorem{theorem}{Theorem}

\makeatletter

\newcommand{\Rom}[1]{\expandafter\@slowromancap\romannumeral #1@}
\makeatother

\title{Modeling and analysis of a novel two-strain dengue epidemics model considering secondary infections with increased mortality}

\author[MAI,IQCB]{Burcu G\"urb\"uz}
\ead{burcu.gurbuz@uni-mainz.de}

\author[ODU]{Ayt\"ul G\"ok\c{c}e}
\ead{aytulgokce@odu.edu.tr}

\author[ESPOL,CFD]{Joseph P\'aez Ch\'avez}
\ead{jpaez@espol.edu.ec}

\author[KOB]{Thomas G\"otz}
\ead{goetz@uni-koblenz.de}

\address[MAI]{Institut f\"ur Mathematik, Johannes Gutenberg-Universit\"at Mainz, Staudingerweg 9, 55099 Mainz, Germany}

\address[ODU]{Ordu University, Faculty of Science and Letters, Department of Mathematics, 52200 Ordu, Turkey}

\address[ESPOL]{Center for Applied Dynamical Systems and Computational Methods (CADSCOM), Faculty of Natural Sciences and Mathematics, Escuela Superior Polit\'ecnica del Litoral, P.O. Box 09-01-5863, Guayaquil, Ecuador}

\address[KOB]{Mathematical Institute, University of Koblenz, D-56070 Koblenz, Germany}

\address[CFD]{Center for Dynamics, Department of Mathematics, TU Dresden, D-01062 Dresden, Germany}

\address[IQCB]{Institute for Quantitative and Computational Biosciences (IQCB), Johannes Gutenberg University Mainz, 55128 Mainz, Germany}

\begin{document}

\begin{abstract}
In this study, we develop and analyze a deterministic two-strain host--vector model for dengue transmission that incorporates key immuno-epidemiological mechanisms, including temporary cross-immunity, antibody-dependent enhancement (ADE), disease-induced mortality during secondary infections, and explicit vector co-infection. The human population is divided into compartments for primary and secondary infections, while the mosquito population includes single- and co-infected classes. ADE is modeled through distinct primary ($\alpha$) and secondary ($\sigma$) transmission rates. Using the next-generation matrix method, we derive the basic reproduction number $R_0$ and establish the local stability of the disease-free equilibrium for $R_0 < 1$. Analytical results show that one-strain endemic equilibria lose stability under ADE conditions ($\sigma > \alpha$), allowing invasion by a heterologous strain. Employing center-manifold theory and numerical continuation (COCO), we demonstrate the occurrence of backward bifurcation, bistability between disease-free and endemic states, and Hopf-induced oscillations. Numerical simulations confirm transitions among disease-free, endemic, and periodic regimes as key parameters vary. The model highlights how ADE, waning cross-immunity, and vector co-infection interact to generate complex dengue dynamics and provides insights useful for designing effective control and vaccination strategies in dengue-endemic regions.
\end{abstract}




\maketitle

\section{Introduction}
\label{Sec:Intro}
Dengue fever is one of the most significant vector-borne diseases worldwide. It causes an estimated 390 million infections and hundreds of thousands of hospitalizations per year \cite{who_dengue}. Secondary infection with a heterologous serotype can lead to severe dengue hemorrhagic fever (DHF) or dengue shock syndrome (DSS). The dengue virus exists as four antigenically distinct serotypes (DENV-1 to DENV-4). While primary infection elicits lifelong immunity to the infecting serotype, cross-immunity to other serotypes wanes rapidly. Additionally, pre-existing antibodies can facilitate infection by a different serotype via antibody-dependent enhancement (ADE). ADE has long been considered a significant risk factor for severe dengue disease during secondary heterotypic infections. There is clinical evidence linking antibody titer ranges to disease severity in children \cite{guzman2010complexity, enwiki:1303851882}.

Mathematical modeling has played a vital role in understanding the epidemiological consequences of multi-strain dynamics, ADE, and transient cross-immunity \cite{janreung2020mathematical}. Early theoretical work demonstrated that simple two-strain models with temporary cross-immunity and immune enhancement can exhibit behaviors such as backward bifurcations, oscillations, and chaos \cite{aguiar2009torus}. For instance, Aguiar et al. (2008) introduced a minimal two-strain SIR-like model showing complex bifurcation structures and chaotic attractors under biologically plausible parameter regimes \cite{aguiar2008epidemiology}. Bianco et al. (2009) analytically explored how weak temporary cross‑immunity may stabilize otherwise chaotic ADE dynamics, and how stronger cross-immunity can induce complex epidemic cycles \cite{bianco2009epidemics}. More recent studies have modified these frameworks to explicitly incorporate vector-host transmission dynamics. Rashkov et al. (2021) developed a compartmental host-vector model incorporating two dengue strains, temporary cross-immunity, and potential secondary infection enhancement \cite{rashkov2021complexity}. This model shows multi-strain equilibrium and bifurcation phenomena consistent with field observations. Ndii et al. (2024) applied a two-serotype host-vector model to quantify the effects of ADE on epidemic persistence. They highlighted the roles of cross-immunity duration and enhancement in shaping outbreak patterns \cite{ndii2024numerical}. Recently, Aguiar et al. (2024) revisited a two-infection SIR-SIR model incorporating temporary immunity and disease enhancement. They uncovered codimension-2 bifurcations that serve as organizing centers for transitions between endemic and multi-strain regimes \cite{aguiar2024bifurcation}. Steindorf et al. (2022) complemented this work by formulating a delay differential model that includes distributed temporary cross-immunity and ADE. They used Lyapunov functionals to explain global stability conditions and periodic outbreak dynamics \cite{steindorf2022cross}.

In this study, we consider a two-strain mathematical model of dengue fever. This deterministic, compartmental model describes the transmission dynamics of two interacting serotypes of the dengue virus within a host-vector framework. The model is based on well-established immunological and epidemiological evidence showing that a dengue infection with one serotype results in lifelong immunity to that serotype, but only temporary cross-immunity to others.

Individuals are at increased risk of developing severe disease upon secondary infection with a different serotype due to ADE \cite{halstead2003neutralization, ferguson1999transmission}. Based on this concept, we present a two-strain mathematical model of dengue fever that explicitly distinguishes between primary and secondary infection classes in humans and co-infection states in mosquitoes. Our framework takes into account temporary cross-immunity, immune enhancement with distinct primary and secondary infection rates \(\alpha\) and \(\sigma\), disease-induced mortality during secondary infection (\(\delta\)), and vector co-infection compartments (\(V_{12}\)). The model allows us to rigorously derive the basic reproduction number,  \(R_0\), analysis of one‑strain and two‑strain endemic equilibrium, and bifurcation behaviors, including backward bifurcation, under ADE conditions (\(\sigma>\alpha\)).

Compared to earlier models, our formulation offers significant improvements that address biological realism and mathematical tractability. First, it integrates explicit vector co-infection states, thereby bridging the gap between purely host-based and classical vector-host models. This allows for a more accurate representation of the full transmission cycle, including sequential infections in vectors. Second, the model mechanistically incorporates ADE by distinguishing between primary and secondary infection pathways with different transmission parameters, denoted by \(\alpha\) and \(\sigma\), respectively. This approach allows for a more realistic exploration of how immune history affects susceptibility and transmission risk during secondary infection events. Third, the model captures waning cross-immunity by introducing intermediate compartments, \(S_1\) and \(S_2\) for partially susceptible individuals who transition from a recovered state back to susceptibility for heterologous strains. This structure better reflects the immunological processes observed in empirical dengue studies. Finally, the model supports rigorous bifurcation analysis and numerical continuation methods. These methods are essential tools for identifying thresholds, bifurcation points, and regions of multi-stability in the system's parameter space. This structured modeling framework aligns with and extends previous literature demonstrating how temporary immunity and ADE interact to shape dengue dynamics. Our primary objective is to clarify the parameter conditions under which phenomena such as multi-strain coexistence, backward bifurcation, and complex epidemic oscillations occur. These insights are critical for developing control strategies and public health interventions in dengue-endemic settings.

The remainder of the paper is organized as follows. Section~\ref{Sec:Model} presents the full mathematical formulation of the two-strain host–vector system. In Section~\ref{Sec:Prelim}, we discuss model positivity, boundedness, and the disease-free equilibrium, and derive $R_0$ via the next-generation matrix method. The one-strain and two-strain endemic equilibria are analyzed in detail, followed by bifurcation and stability analysis. Section~\ref{Sec:NumSim} presents numerical continuation results illustrating transitions between equilibria and oscillatory dynamics. Finally, Section~\ref{Sec:Conc} concludes with implications for dengue control strategies and future model extensions.

\section{The mathematical model}
\label{Sec:Model}

In this model, we assume total human and vector populations of size  $N$ and $M$. The human population is divided into several compartments: susceptibles, $S$ individuals infected with strain 1, $I_1$ or strain 2, $I_2$, individuals recovered from strain 1, $R_1$, or strain 2, $R_2$, individuals susceptible after recovery, $S_1$, $S_2$, individuals reinfected with strain 2 after strain 1, $I_{12}$, or strain 1 after strain 2, $I_{21}$, and fully immune individuals, $R$. The mosquito population includes susceptible vectors, $U$, and vectors infected with strain 1, $V_1$, strain 2, $V_2$, or co-infected, $V_{12}$. The model represents primary infection with terms involving the transmission rate, $\alpha$, from infectious vectors to susceptible humans. Those who recover from a primary infection  $I_1$ or $I_2$ transition to the temporary immune classes $R_1$ and $R_2$. From there, they move into the partially susceptible compartments $S_1$ and $S_2$ at rate $\nu$. These individuals are then at risk of a secondary infection, which is represented by the compartments $I_{12}$ and $I_{21}$ with a transmission rate of $\sigma$, which may differ from $\alpha$ to reflect ADE effects. The mosquito dynamics are modeled through interactions with all infected human classes. Susceptible vectors, $U$, acquire infection at a rate proportional to the prevalence of human infection and transition into single or co-infected states. The co-infection compartment $V_{12}$ represents vectors exposed to both strains through sequential biting events, an abstraction of superinfection and viral interference \cite{adams2006cross}.

\begin{figure}[H]
\begin{tikzpicture}[
  node distance=1.7cm and 1.5cm,
  compartment/.style={circle, draw, fill=cyan!30, minimum size=1.0cm, font=\small},
  arrow/.style={-{Latex[length=2mm]}, thick},
  redarrow/.style={arrow, red},
  greenarrow/.style={arrow, green!60!black},
  bluearrow/.style={arrow, blue},
  dashedarrow/.style={arrow, dashed}
]

\node[compartment] (S) {S};
\node[compartment, right=of S, yshift=1.25cm] (I1) {$I_1$};
\node[compartment, right=of I1] (R1) {$R_1$};
\node[compartment, right=of R1] (S1) {$S_1$};
\node[compartment, right=of S1] (I12) {$I_{12}$};
\node[compartment, right=of I12, yshift=-1.25cm] (R) {R};

\node[compartment, below=of I1] (I2) {$I_2$};
\node[compartment, right=of I2] (R2) {$R_2$};
\node[compartment, right=of R2] (S2) {$S_2$};
\node[compartment, right=of S2] (I21) {$I_{21}$};

\draw[redarrow] (S) -- (I1) node[midway, below] {\hspace*{2em}\scriptsize $\alpha(V_1 + V_{12})$};
\draw[redarrow] (S) -- (I2) node[midway, below left=-1pt] {\scriptsize $\alpha(V_2 + V_{12})$};

\draw[greenarrow] (I1) -- (R1) node[midway, above] {\scriptsize $\gamma$};
\draw[greenarrow] (I2) -- (R2) node[midway, below] {\scriptsize $\gamma$};

\draw[bluearrow] (R1) -- (S1) node[midway, above] {\scriptsize $\nu$};
\draw[bluearrow] (R2) -- (S2) node[midway, below] {\scriptsize $\nu$};

\draw[redarrow] (S1) -- (I12) node[midway, above] {\scriptsize $\sigma(V_2 + V_{12})$};
\draw[redarrow] (S2) -- (I21) node[midway, below] {\scriptsize $\sigma(V_1 + V_{12})$};

\draw[greenarrow] (I12) -- (R) node[midway, above] {\scriptsize $\gamma$};
\draw[greenarrow] (I21) -- (R) node[midway, below right] {\scriptsize $\gamma$};

\draw[dashedarrow] (I12) -- ++(0,1) node[right] {\scriptsize $\delta$};
\draw[dashedarrow] (I21) -- ++(0,-1) node[right] {\scriptsize $\delta$};

\foreach \x in {S,I1,I2,R1,R2,S1,S2,I12,I21,R}
  \draw[arrow] (\x) -- ++(0,1.1) node[above] {\scriptsize $\mu$};

\draw[arrow] ++(-1.3,0) node[left] {\scriptsize $\Lambda_N$} -- (S);

  node distance=2cm and 2cm,
  virus/.style={circle, draw=black, fill=gray!30, minimum size=1.4cm, font=\large},
  arrow/.style={-{Latex[length=2mm]}, thick},
  redarrow/.style={arrow, red, line width=1pt},
  blackarrow/.style={arrow, black},
  ]
\end{tikzpicture}
\centering
\begin{tikzpicture}[
  node distance=1.5cm and 1.0cm,
  virus/.style={circle, draw=black, fill=gray!30, minimum size=1.0cm, font=\large},
  arrow/.style={-{Latex[length=2mm]}, thick},
  redarrow/.style={arrow, red, line width=1pt},
  blackarrow/.style={arrow, black},
  ]
\node[virus] (U) {$U$};
\node[virus, above right=of U] (V1) {$V_1$};
\node[virus, below right=of U] (V2) {$V_2$};
\node[virus, right=of V1, yshift=-2.25cm] (V12) {$V_{12}$};

\draw[blackarrow] (U) -- ++(135:1.4) node[midway, left] {\scriptsize $\kappa$};
\draw[blackarrow] (U) ++(0,-1.2) -- (U) node[midway, right] {\scriptsize $\Lambda_M$};

\draw[redarrow] (U) -- (V1) node[midway, above left=-2pt] {\scriptsize $\beta(I_1 + I_{21})$};
\draw[redarrow] (U) -- (V2) node[midway, below left=-2pt] {\scriptsize $\beta(I_2 + I_{12})$};

\draw[redarrow] (V1) -- (V12) node[midway, above right=-2pt] {\scriptsize $\beta(I_2 + I_{12})$};
\draw[redarrow] (V2) -- (V12) node[midway, below right=-2pt] {\scriptsize $\beta(I_1 + I_{21})$};

\draw[blackarrow] (V1) -- ++(45:1.5) node[midway, right] {\scriptsize $\kappa$};
\draw[blackarrow] (V2) -- ++(-45:1.5) node[midway, right] {\scriptsize $\kappa$};
\draw[blackarrow] (V12) -- ++(0:1.5) node[midway, above] {\scriptsize $\kappa$};
\end{tikzpicture}
\caption{
Schematic representation of the host-vector transmission dynamics of a two-strain pathogen. The human population \( N \) is divided into compartments based on infection status: susceptible (\( S \)), infected with strain \( i \) (\( I_i \)), recovered from strain \( i \) (\( R_i \)), partially susceptible after infection with strain \( i \) (\( S_i \)), reinfected with strain \( j \) after recovering from strain \( i \) (\( I_{ij} \)), and fully recovered (\( R \)) where $i,j \in 1,2$. The vector population includes susceptible vectors (\( U \)) and vectors infected with strain 1, strain 2, or both (\( V_1, V_2, V_{12} \)). Arrows indicate transitions due to infection, recovery, waning immunity, and disease-induced mortality.
}
\label{Fig:Schem}
\end{figure}
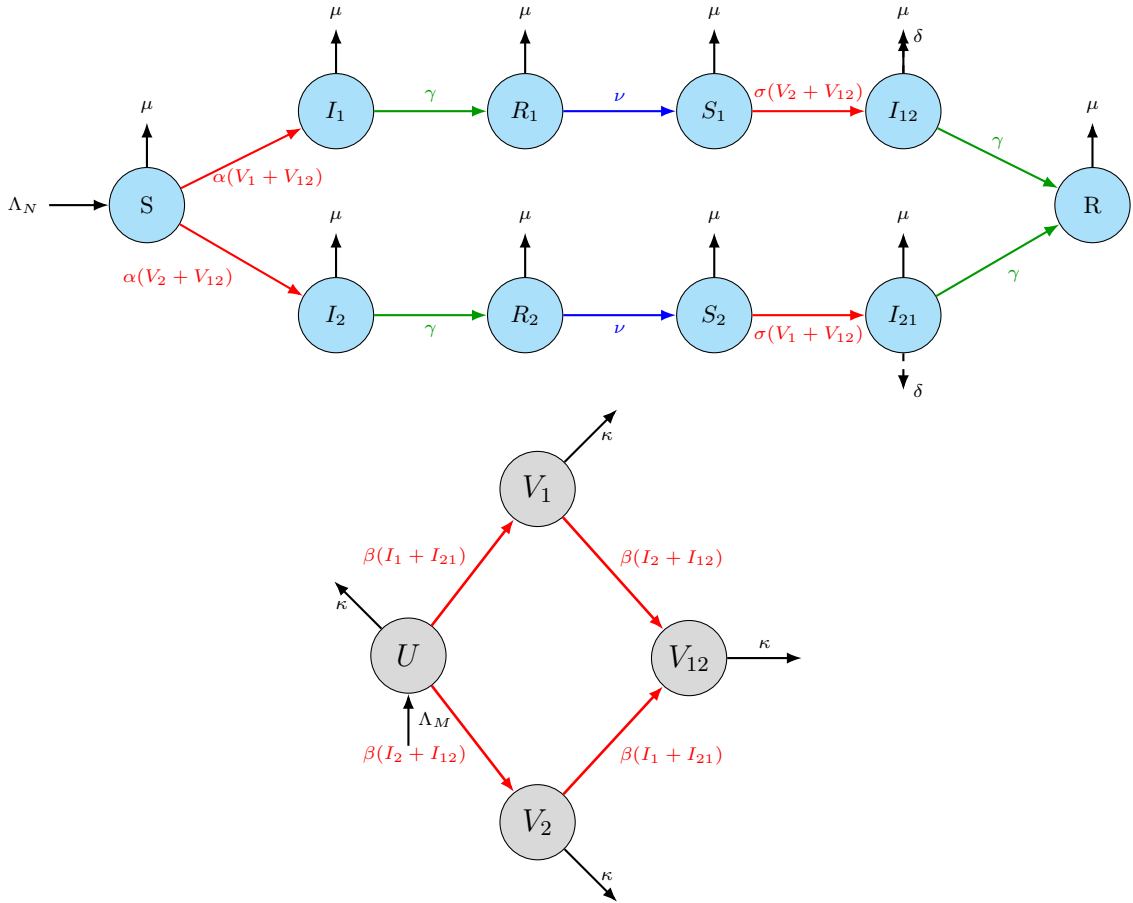


Disease-induced mortality during secondary infection is included via the parameter $\delta$, which can be set to zero for analytic tractability. The vector infection rate $\beta$ governs transmission from infected humans to mosquitoes. The schematic representation of the host-vector transmission dynamics of a two-strain pathogen is shown in Figure \ref{Fig:Schem}.\\

The model is given by
\begin{subequations}
\label{Eq:Model}
\begin{align}
    \Dot{S} &= \Lambda_N - \frac{\alpha}{M} S (V_{1}+V_{2}+2 V_{12})-\mu S, \label{Eq:S}\\
    \Dot{I}_1 &= \frac{\alpha}{M} S (V_{1}+V_{12}) - (\gamma+\mu) I_1, \label{Eq:I1}\\
    \Dot{I}_2 &= \frac{\alpha}{M} S (V_{2}+V_{12}) - (\gamma+\mu) I_2,  \label{Eq:I2}\\
    \Dot{R}_1 &= \gamma I_1 - (\nu+\mu) R_1, \label{Eq:R1}\\
    \Dot{R}_2 &= \gamma I_2 - (\nu+\mu) R_2, \label{Eq:R2}\\
    \Dot{S}_1 &= \nu R_1 - \frac{\sigma}{M} S_1 (V_{2}+V_{12}) -\mu S_1, \label{Eq:S1}\\
    \Dot{S}_2 &= \nu R_2 - \frac{\sigma}{M} S_2 (V_{1}+V_{12})-\mu S_2, \label{Eq:S2}\\
    \Dot{I}_{12} &= \frac{\sigma}{M} S_1 (V_2+V_{12}) - (\gamma+\mu) I_{12} - \delta I_{12}, \label{Eq:I12}\\
    \Dot{I}_{21} &= \frac{\sigma}{M} S_2 (V_1+V_{12}) - (\gamma+\mu) I_{21} - \delta I_{21}, \label{Eq:I21}\\
    \Dot{R} &= \gamma (I_{12}+I_{21}) -\mu R, \label{Eq:R}
\intertext{and}
    \Dot{U} &= \Lambda_M -\frac{\beta}{N} U (I_1+I_2+I_{12}+I_{21}) - \kappa U, \label{Eq:U}\\
    \Dot{V}_1 &= \frac{\beta}{N} U (I_1+I_{21}) -\frac{\beta}{N} V_1(I_2+I_{12}) - \kappa V_1, \label{Eq:V1}\\
    \Dot{V}_2 &= \frac{\beta}{N} U (I_2+I_{12}) -\frac{\beta}{N} V_2(I_1+I_{21}) - \kappa V_2, \label{Eq:V2}\\
    \Dot{V}_{12} &= \frac{\beta}{N} V_1 (I_2 + I_{12}) + \frac{\beta}{N} V_2 (I_1 +I_{21}) - \kappa V_{12}. \label{Eq:V12}
\end{align}
\end{subequations}

The model parameters are informed by epidemiological data. The average human lifetime is taken as $\mu^{-1} = 12 \times 65$ month unit, the recovery period is $\gamma^{-1} = 0.5$ month unit, and the cross-immunity period is $\nu^{-1} = 9$ month unit \cite{ferguson1999transmission}. The average vector lifespan is $\kappa^{-1} = 1$–1.5 month unit. A possible death induced mortality rate is called $\delta$. For the beginning, we may assume $\delta=0$ and hence a \textbf{constant} human population. The human infection rates $\alpha, \sigma$ for the primary and secondary infection as well as the vector infection rate $\beta$ are less clear. In particular, the ADE effect is modeled by varying $\sigma$ relative to $\alpha$. Three biological scenarios are typically examined:
\begin{itemize}
    \item $\sigma < \alpha$: secondary infections are less likely than primary infections, implying partial protection.
    \item $\sigma = \alpha$: equal risk of primary and secondary infections.
    \item $\sigma > \alpha$: ADE enhances susceptibility to secondary infection.
\end{itemize}
The model variables and parameters are shown in Table \ref{Tab:VarPar}.

\begin{table}[H]
    \centering
    \caption{Description of the variables and parameters of the  model \eqref{Eq:Model}.}
    \begin{tabularx}{\textwidth}{c | l}
    \hline  \hline
		Variable & Description\\
		\hline  \hline
$N$ & Total human population\\
$M$ & Total vector population\\
$S$ & Fully susceptible (never infected)\\
$I_1$, $I_2$ & Primary infected with strain 1 or 2, respectively\\
$R_1$, $R_2$ &  Recovered from strain 1 or 2, with temporary immunity\\
$S_1$, $S_2$ & Partially susceptible after immunity wanes from infection with strain 1 or 2\\
$I_{12}$, $I_{21}$ & Secondary infected, previously recovered from strain 1 (or 2), now infected by strain 2 (or 1)\\
$R$ &  Fully recovered after two infections\\
$U$ &  Uninfected mosquitoes\\
$V_1$ &  Infected with strain 1\\
$V_2$ &  Infected with strain 2\\
$V_{12}$ &  Co-infected with both strains\\ 

		\hline \hline
		Parameter & Description\\
		\hline \hline

$\Lambda_N$, $\Lambda_M$ & Recruitment of humans/vectors \\
$\mu$, $\kappa$ & Natural death rates of humans and vectors \\

$\alpha$ & Primary infection rate from infected vectors\\
$\sigma$ & Secondary infection rate (may differ from
 $\alpha$)\\
$\gamma$ & Human recovery rate\\

$\nu$ & Rate of loss of temporary immunity\\
$\delta$ & Disease-induced mortality in secondary infections\\

$\beta$ & Transmission from infected humans to vectors\\
		\hline \hline
    \end{tabularx}
	\label{Tab:VarPar}
\end{table}

This model builds on earlier multi-strain and host-vector frameworks \cite{kooi2014analysis, kooi2013bifurcation, rashkov2021complexity}, while incorporating enhanced biological realism through explicit co-infection classes and multiple infection stages. These models are essential for understanding the nonlinear dynamics of dengue outbreaks, determining the role of ADE in epidemic amplification, and developing public health strategies, including vector control and vaccination.

In particular, the present formulation explicitly distinguishes between primary and secondary human infections and allows mosquitoes to exist in single- or double-infected states. The inclusion of the co-infected vector class $V_{12}$ enables the description of sequential acquisition of different serotypes by the same vector, which is an essential feature for understanding viral competition and coexistence within the mosquito population. This mechanism captures both the potential for viral interference and the possibility of enhanced transmission when mosquitoes harbor multiple viral strains simultaneously \cite{adams2006cross}. Such vector–host coupling provides a more accurate framework than purely host-based models and allows the system to reproduce the epidemiological complexity observed in multi-strain dengue outbreaks.

From an immuno-epidemiological perspective, the compartments $S_1$ and $S_2$ play a critical role in describing individuals who have lost temporary cross-immunity and have returned to a partially susceptible state. The transition from $R_i$ to $S_i$ at rate $\nu$ reflects the finite duration of heterologous protection following recovery from a primary infection. This process introduces a biologically relevant delay between infections, enabling the model to capture realistic inter-epidemic intervals observed in dengue-endemic regions. The secondary infection process, parameterized by the transmission rate $\sigma$, accounts for ADE effects by modifying susceptibility to reinfection. The ratio $\sigma / \alpha$ thus quantifies the degree of enhancement or protection during heterologous exposure, where $\sigma > \alpha$ corresponds to enhanced transmission due to ADE, and $\sigma < \alpha$ corresponds to partial cross-protection.

The inclusion of disease-induced mortality, represented by $\delta$, allows the model to account for severe outcomes of secondary infections such as DHF and DSS. This term introduces an additional loss pathway for infected individuals and provides a means to assess how mortality alters disease persistence and endemic equilibria. Although $\delta$ is typically small in magnitude, its impact on model dynamics can be nontrivial, especially when combined with ADE-driven secondary infection amplification.

Overall, the system \eqref{Eq:Model} provides a unified mechanistic description of multi-strain dengue dynamics that links immunological processes at the individual level (temporary cross-immunity, ADE, secondary mortality) with transmission mechanisms at the population level (vector co-infection, human reinfection). This structure not only facilitates rigorous mathematical analysis of equilibria and bifurcation behavior but also supports biologically meaningful interpretation of parameters. The model thus serves as a comprehensive foundation for exploring threshold phenomena, coexistence of dengue strains, and oscillatory epidemic patterns under varying epidemiological and environmental conditions.

\section{Analytical Results} 
\label{Sec:Prelim}
For any initial values at a given time, the right hand side of the system \eqref{Eq:S}-\eqref{Eq:V12} has a non-negative solution for $t>0$. Let $N(t)=S(t)+I_1(t)+I_2(t)+R_1(t)+R_2(t)+S_1(t)+S_2(t)+I_{12}(t)+I_{21}(t)+R(t)$ and $M(t)=U(t)+V_1(t)+V_2(t)+V_{12}(t)$. Then
\begin{align}
    \frac{d N}{d t} &= \Lambda_N-\mu N-\delta(I_{12}+I_{21})\leq \Lambda_N-\mu N,\nonumber\\
    \frac{d M}{d t} &= \Lambda_M-\kappa M.\nonumber
\end{align}
This implies that $N$ and $M$ remains bounded on any finite time interval. Using differential inequality for $N$ and $M$ leads to
\begin{equation}
    \limsup_{t\to\infty} N(t)\leq \frac{\Lambda_N}{\mu},\quad \textrm{and}\quad \limsup_{t\to\infty} M(t) \leq \frac{\Lambda_M}{\kappa}.\nonumber
\end{equation}

\subsection{Disease-free Equilibrium (DFE)}
\label{Subsec:DFE}
The trivial equilibrium of the model \eqref{Eq:S}-\eqref{Eq:V12} is $E_0=(S^*,0,0,0,0,0,0,0,0,0,U^*,0,0,0)$ with $S^*=\Lambda_N/\mu$ and $U^*=\Lambda_M/\kappa.$
Linearization around the disease-free steady state leads to the Jacobian matrix
where Jacobian can be re-expressed as
$$\mathcal{J}=[\mathcal{M}_{{14}\times7} \;\;\mathcal{N}_{{14}\times7}]$$
with
$$
\mathcal{M}=\begin{bmatrix}
-\mu & 0 & 0 & 0 & 0 & 0 & 0 \\
 0& -\bar{\gamma} & 0 & 0 & 0 & 0 & 0 \\
 0 & 0 & -\bar{\gamma} & 0 & 0 & 0 & 0 \\
 0 & \gamma & 0 & -\bar{\nu} & 0 & 0 & 0 \\
 0 & 0 & \gamma & 0 & -\bar{\nu} & 0 & 0 \\
 0 & 0 & 0 & \nu & 0 & -\mu & 0 \\
 0 & 0 & 0 & 0 & \nu & 0 & -\mu \\
 0 & 0 & 0 & 0 & 0 & 0 & 0 \\
 0 & 0 & 0 & 0 & 0 & 0 & 0\\
 0 & 0 & 0 & 0 & 0 & 0 & 0 \\
  0 & -\frac{\beta U^*}{S^*} & -\frac{\beta U^*}{S^*} & 0 & 0 & 0 & 0 \\
   0 & \frac{\beta U^*}{S^*} & 0 & 0 & 0 & 0 & 0 \\
0 & 0 & \frac{\beta U^*}{S^*} & 0 & 0 & 0 & 0  \\
0 & 0  & 0  & 0 & 0 & 0 & 0 \\
\end{bmatrix}
$$
and
$$
\mathcal{N}=\begin{bmatrix}
0 & 0 & 0 & 0 & -\frac{\alpha S^*}{U^*} & -\frac{\alpha S^*}{U^*} & -\frac{2 \alpha S^*}{U^*} \\
0 & 0 & 0 & 0 & \frac{\alpha S^*}{U^*} & 0 & \frac{\alpha S^*}{U^*}\\
 0 & 0 & 0 & 0 & 0 & \frac{\alpha S^*}{U^*} & \frac{\alpha S^*}{U^*}\\
 0 & 0 & 0 & 0 & 0 & 0 & 0\\
 0 & 0 & 0 & 0 & 0 & 0 & 0\\
0 & 0 & 0 & 0 & 0 & 0 & 0\\
 0 & 0 & 0 & 0 & 0 & 0 & 0\\
 -\bar{\delta} & 0 & 0 & 0 & 0 & 0 & 0\\
 0 & -\bar{\delta} & 0 & 0 & 0 & 0 & 0\\
  \gamma & \gamma & -\mu & 0 & 0 & 0 & 0\\
 -\frac{\beta U^*}{S^*} & -\frac{\beta U^*}{S^*} & 0 & -\kappa & 0 & 0 & 0\\
 0 & \frac{\beta U^*}{S^*} & 0 & 0 & -\kappa  & 0 & 0\\
  \frac{\beta U^*}{S^*}&0  & 0 & 0 & 0 & -\kappa & 0 \\
  0 & 0  & 0 & 0 & 0 & 0 & -\kappa\\
\end{bmatrix}
$$
where $\bar{\nu}=\nu+\mu$, $\bar{\delta}=\gamma+\mu+\delta$.
Here the characteristic polynomial $\textrm{det}(J|_{E_0}-\Sigma I)=0$ results in eigenvalues corresponding to the disease-free steady state:$-\mu, -\kappa, -(\mu+\nu), -(\mu+\gamma+\delta)$ (counted with multiplicity) and
\begin{equation}
    \Sigma_{11:14}= -\frac{1}{2} \left(\gamma+\kappa+\mu)\pm \left((\gamma+\kappa+\mu)^2-4(\gamma \kappa+\mu \kappa-\alpha \beta)\right)^{1/2}\right), \nonumber
\end{equation}
which are the roots of polynomial
\begin{equation}
    \Sigma^2+(\gamma+\kappa+\mu)\Sigma+ \gamma\kappa+\kappa \mu-\alpha \beta=0.\label{Eq:DFE_Poly}
\end{equation}
The linearization around a generic steady state is given in \ref{Sec:Appendix1}. The polynomial given in Eq. \eqref{Eq:DFE_Poly} has a positive real root if the free term is negative, namely
\begin{equation}
    \frac{\alpha \beta}{\kappa(\gamma+\mu)}>1.\nonumber
\end{equation}
This is potentially associated with the reproduction number

\subsection{Basic reproduction number}
\label{Subsec:R0}
In this section, following the ideas given in \cite{van2002reproduction}, the basic reproduction number for \eqref{Eq:S}–\eqref{Eq:V12} is derived using the next generation matrix method. For disease free state, the matrix $F$ associated with new
infections and the matrix $V$ representing the remaining expressions are given by
$$F=
\begin{bmatrix}
0& 0 & 0 & 0 & \frac{\alpha S^*}{M^*} & 0 & \frac{\alpha S^*}{M^*}  \\
0& 0 & 0 & 0 & 0 & \frac{\alpha S^*}{M^*}  &\frac{\alpha S^*}{M^*}  \\
0& 0 & 0 & 0 & 0& 0&0 \\
0& 0 & 0 & 0 & 0& 0&0 \\
\frac{\beta U^*}{N^*}& 0 & 0 & \frac{\beta U^*}{N^*} & 0&0&0 \\
0& \frac{\beta U^*}{N^*} & \frac{\beta U^*}{N^*} & 0&0&0&0 \\
0& 0 & 0 & 0 & 0& 0&0
\end{bmatrix}
$$
$$V=-\mathcal{V}=
\begin{bmatrix}
\gamma+\mu& 0 & 0 & 0 & 0& 0&0 \\
0& \gamma+\mu & 0 & 0 & 0&0&0 \\
0& 0 & \gamma+\mu+\delta & 0 & 0& 0&0 \\
0& 0 & 0 & \gamma+\mu+\delta & 0& 0&0 \\
0& 0 & 0 & 0 & \kappa&0&0 \\
0&0& 0 & 0&0&\kappa&0 \\
0& 0 & 0 & 0 & 0& 0&\kappa
\end{bmatrix}.
$$
Since the inverse of a diagonal matrix is given by replacing the main diagonal components of the matrix with their reciprocals, we have
$$FV^{-1}=
\begin{bmatrix}
0& 0 & 0 & 0 & \frac{\alpha S^*}{\kappa M^*}& 0 & \frac{\alpha S^*}{\kappa M^*} \\
0& 0 & 0 & 0 & 0& \frac{\alpha S^*}{\kappa M^*} & \frac{\alpha S^*}{\kappa M^*}\\
0& 0 & 0 & 0 & 0& 0&0 \\
0& 0 & 0 & 0 & 0& 0&0 \\
\frac{\beta U^*}{ (\gamma+\mu)N^* }& 0 & 0 & \frac{\beta U^*}{ (\gamma+\mu+\delta)N^*} & 0&0&0 \\
0&\frac{\beta U^*}{ (\gamma+\mu)N^* } &   \frac{\beta U^*}{ (\gamma+\mu+\delta)N^*} & 0&0&0&0 \\
0& 0 & 0 & 0 & 0& 0&0
\end{bmatrix}.
$$
The reproduction number $R_0$ is defined  to be the spectral radius of the matrix $FV^{-1}$. The characteristic equation of this matrix is given by $\textrm{det}(FV^{-1}-PI)=0$, leading to
\begin{equation}
(P^4-2A_1A_2P^2+A_1^2A_2^2)P^3=0,\label{Eq:Eigen_R0}
\end{equation}
where $A_1= \frac{\alpha S^*}{\kappa M^*}$ and $A_2=\frac{\beta U^*}{  (\gamma+\mu) N^*}$. The nonzero roots of Eq. \eqref{Eq:Eigen_R0} are $P_{1,2}=\pm \sqrt{A_1A_2}$. Hence the basic reproduction number is associated with the dominant eigenvalue $P_1$. Replacing $N^*=S^*=\frac{\Lambda_N}{\mu}$ and $M^*=U^*=\frac{\Lambda_M}{\kappa}$ leads to
\begin{equation}
    R_0=\sqrt{A_1A_2}=\left(\frac{\alpha \beta }{\kappa (\gamma+\mu)}\right)^{1/2}.
\end{equation}

\begin{theorem}
The basic reproduction number \( \mathcal{R}_0 \) of the system \eqref{Eq:Model} is given by
\[
\mathcal{R}_0 = \left( \frac{\alpha \beta}{\kappa (\gamma + \mu)} \right)^{1/2}.
\]
The disease-free equilibrium is locally asymptotically stable if \( \mathcal{R}_0 < 1 \), and unstable if \( \mathcal{R}_0 > 1 \) \cite{van2002reproduction}.
\end{theorem}


\subsection{One-strain equilibrium}
Let $E_1$ denote the equilibrium where the second dengue strain is absent, namely

$$E_1=\left(S_{01}^*,I_1^*,0,R_1^*,0,S_1^*,0,0,0,0,U_0^*,V_1^*,0,0\right).$$
In this context, the model \eqref{Eq:S}-\eqref{Eq:V12} leads to
\begin{align}
    I_1^*&=\frac{\alpha}{(\gamma+\mu)M_{*}} S^* V_1^*,\label{Eq:I1s}\\
    R_1^*&=\frac{\gamma}{(\nu+\mu)}I_1^*=\frac{\gamma \alpha}{(\gamma+\mu)(\nu+\mu)M_{*}} S^* V_1^*,\label{Eq:R1s}\\
    S_1^*&=\frac{\nu}{\mu}R_1^* = \frac{\nu\gamma\alpha}{\mu(\gamma+\mu)(\nu+\mu)M_{*}} S^* V_1^*,\label{Eq:S1s}
\end{align}
with $N_{*}=S_{01}^*+I_1^*+R_1^*+S_1^*$ and $M_{*}=U_0^*+V_1^*$. Letting $\mathcal{X}_1=S_{01}^*$ and $\mathcal{X}_2=S_{01}^* V_1^*$, we obtain
\begin{align}
    N_{*}=&\mathcal{X}_1 +\frac{\alpha}{(\gamma+\mu)M_{*}}\left(1+\frac{\gamma}{(\nu+\mu)}+\frac{\nu\gamma}{\mu(\nu+\mu)}\right)\mathcal{X}_2, \notag\\
    =&\mathcal{X}_1+\frac{\alpha}{\mu M_{*}}\mathcal{X}_2.\label{Eq:X1X2}
\end{align}
Using Eqs. \eqref{Eq:V1} and \eqref{Eq:I1s}, one obtains
\begin{equation}
    \kappa=\frac{\beta}{N_{*}} \left(\frac{\alpha}{(\gamma+\mu)} \mathcal{X}_1-\frac{\alpha}{(\gamma+\mu)M_{*}}\mathcal{X}_2\right). \label{Eq:X1X22}
\end{equation}
Solving \eqref{Eq:X1X2} and \eqref{Eq:X1X22} leads to
\begin{align}
    \mathcal{X}_1&= \frac{N_{*}}{\alpha+\mu} \left(\mu+ \frac{\kappa(\gamma+\mu)}{\beta}\right),\\
    \mathcal{X}_2&=\frac{\mu M_{*}N_{*}}{\alpha+\mu}\left(1-\frac{\kappa(\gamma+\mu)}{\alpha\beta}\right).\label{Eq:X2}
\end{align}
Therefore considering $R_0=\left(\frac{\alpha\beta}{\kappa(\gamma+\mu)}\right)^{1/2}$ and using Eqs. \eqref{Eq:I1s}-\eqref{Eq:S1s}, the non-zero components for human population
\begin{align}
    S_{01}^* &=\frac{N_{*}}{\alpha+\mu} \left(\mu+\frac{\kappa(\gamma+\mu)}{\beta}\right)
        = \frac{\alpha N_*}{\alpha+\mu} \left(\frac{\mu}{\alpha}+\frac{1}{R_0^2}\right),\label{Eq:S01}\\
    I_1^* &= \frac{\alpha \mu N_{*}}{(\gamma+\mu)(\alpha+\mu)}\left(1-\frac{1}{R_0^2}\right),\label{Eq:I1ss}\\
    R_1^*&= \frac{\alpha\mu\gamma N_{*}}{(\nu+\mu)(\gamma+\mu)(\alpha+\mu)}\left(1-\frac{1}{R_0^2}\right),\\
    S_1^*&=\frac{\alpha\gamma\nu N_{*}}{(\nu+\mu)(\gamma+\mu)(\alpha+\mu)}\left(1-\frac{1}{R_0^2}\right).
\end{align}
On the other hand, using Eqs. \eqref{Eq:U} and \eqref{Eq:I1ss} one obtains
\begin{align}
    U_0^* =\frac{\lambda_M }{\kappa+\varphi},\quad \varphi=\frac{\alpha \beta \mu}{\bar{\gamma}\bar{\alpha}}\left(1-\frac{1}{R_0^2}\right), \label{Eq:U0s}
\end{align}
where $\bar{\gamma}=\gamma+\mu$ and $\bar{\alpha}=\alpha+\mu$. Besides, using Eq. \eqref{Eq:I1ss} in \eqref{Eq:V1} leads to
\begin{align}
    V_1^*=\frac{\Lambda_M \varphi}{\kappa (\kappa+\varphi)}\;
\end{align}
and thus
\begin{align}
    S_{01}^*=\frac{\Lambda_N (\kappa+\varphi)}{\left(\mu(\kappa+\varphi)+\alpha \varphi \right)},\; I_2^*=\frac{\Lambda_N \alpha \varphi}{\bar{\gamma}\left(\mu(\kappa+\varphi)+\alpha \varphi \right)},\; R_2^*=\frac{\gamma}{\bar{\nu}}I_2^*,\; S_2^*=\frac{\nu \gamma}{\bar{\nu}\mu}I_2^*,
\end{align}
with $\bar{\nu}=\nu+\mu.$

%
Similar calculations can be performed for the equilibrium where the second dengue strain is absent ($E_2$). Note that the components of one-strain equilibrium $E_1$ and $E_2$ are positive if and only if $R_0>1$.

\subsection{Two-strain equilibrium}
Now we consider that both virus strains are present in the population. Denoting any arbitrary two-strain equilibrium point of the model \eqref{Eq:S}-\eqref{Eq:V12} as
\[E= \left(S^{**}, I_1^{**}, I_2^{**}, R_1^{**}, R_2^{**}, S_1^{**}, S_2^{**}, I_{12}^{**}, I_{21}^{**}, R^{**},U^{**}, V_1^{**}, V_2^{**}, V_{12}^{**}\right)\]
and considering
\[y_k=\frac{\alpha}{M^{**}}V_k,\quad k=1,2,12,\]
then we obtain
\[S^{**}=\frac{\Lambda_N}{\mu+y_1+y_2+y_{12}}.\]
Assuming that a two-strain equilibrium emerges with the assumption $y_1=y_2=y$ and $x=y+y_{12}$, leading to
\begin{align*}
    I_1^{**} &= I_2^{**}=\frac{2\Lambda_N x}{\bar{\gamma}(\mu+2 x)},\\
    R_1^{**} &= R_2^{**}=\frac{2\gamma \Lambda_N x}{\bar{\nu}\bar{\gamma} (\mu+2 x)},\\
    S_1^{**} &= S_2^{**}=\frac{\nu \gamma \Lambda_N x}{\bar{\nu}\bar{\gamma} (\mu+2x)(\mu+px)},\\
    I_{12}^{**} &= I_{21}^{**}=\frac{p\nu\gamma \Lambda_N x^2}{\bar{\nu} \bar{\gamma} \bar{\delta}(\mu+2x)(\mu+px) },\\
    R^{**} &= \frac{2 p \nu \gamma^2 \Lambda_N x^2}{\mu \bar{\nu} \bar{\gamma} \bar{\delta}(\mu+2x)(\mu+px) },\quad U^{**} = \frac{\Lambda_M}{\kappa+2\frac{\beta}{N} f(x)}, \\
    V_1^{**}&=V_2^{**}= \frac{\frac{\beta}{N} U^{**} f(x)}{\kappa+\frac{\beta}{N^{**}} f(x)}, \quad V_{12}=\frac{2\beta}{\kappa N^{**}} V_1^{**} f(x),
\end{align*}
where
$$f(x)= \frac{\Lambda_N x}{\bar{\gamma} (\mu+2 x)} \left[2+\frac{p q x}{\mu+p x}\right]$$
with $p=\frac{\sigma}{\alpha}$ and $q=\frac{\nu \gamma}{\bar{\nu} \bar{\delta}}.$

Using Eqs. \eqref{Eq:V1} and \eqref{Eq:V12} leads to
\begin{align}
    \frac{\beta}{N^{**}} f(x) \left[\frac{\Lambda_M}{\kappa+\frac{2\beta}{N^{**}} f(x)}-\frac{M^{**}}{\alpha} y\right] - \frac{\kappa M^{**}}{\alpha}y&=0, \label{Eq:D1}\\
    \frac{2\beta}{N^{**}} y f(x)-\kappa(x-y)&=0.\label{Eq:D2}
\end{align}
Solving Eqs. \eqref{Eq:D1} and \eqref{Eq:D2} for $x$ and $y$ variables, it is possible to obtain $y$ and $y_{12}=x-y$. Therefore, a two-strain equilibrium (with the assumption of $y_1$=$y_2$) can be obtained as a function of $M^{**}$ and $N^{**}$. Note that there may exist various equilibrium where $y_1\neq y_2$, yet the complexity of the model does not allow further analytical investigation.

\subsection{Backward bifurcation analysis}
Using center manifold theory, the backward bifurcation analysis can be carried out \cite{castillo2004dynamical, carr2012applications}. Following the ideas given in Theorem 4.1 in \cite{castillo2004dynamical}, the following change of variables is considered:
$$(S, I_1, I_2, R_1, R_2, S_1, S_2, I_{12}, I_{21}, R, U, V_1, V_2, V_{12})=(x_1,x_2,\dots,x_{14})$$
Considering that $\mathcal{X}=(x_1,x_2,\dots,x_{14})^T$, the model can be rewritten as
$$\frac{d\mathcal{X}}{d t}=F(\mathcal{X},p)$$
where $F=(F_1,F_2,\dots,F_{14})$ and $p$ denotes the parameter set. Now, consider the case for which $R_0=1$ and choose $\alpha$ as bifurcation parameter. Solving the parameter $\alpha$ from $R_0=1$ leads to
\begin{equation*}
    \alpha^*=\frac{\kappa(\gamma+\mu)}{\beta}.
\end{equation*}
Using center manifold theory, the dynamics of the system can be analysed in the vicinity of $\alpha=\alpha^*$. In this context, it is necessary to compute the right and left eigen vectors of Jacobian evaluated at disease free equilibrium. Denote the right eigenvector of $J|_{E_0}$ as $\textbf{w}=(\omega_1,\omega_2,\dots,\omega_{14})^T$ where
\begin{equation}\label{E:rvec}
\begin{aligned}
    \omega_1&= -\frac{\bar{\gamma}}{\mu} (\omega_2+\omega_3), \quad \omega_4=\frac{\gamma }{\bar{\nu}}\omega_2, \quad \omega_5=\frac{\gamma }{\bar{\nu}}\omega_3, \quad \omega_6=\frac{\nu}{\mu}\omega_4=\frac{\nu\gamma }{\mu\bar{\nu}}\omega_2,\\
    \omega_7&=\frac{\nu}{\mu}\omega_5=\frac{\nu\gamma }{\mu\bar{\nu}}\omega_3,\quad \omega_8=0,\;\;\omega_9=0,\;\;\omega_{10}=0,\quad \omega_{11}=-\frac{\bar{\gamma}U^*}{\alpha^* S^*} (\omega_2+\omega_3),\\
    \omega_{12}&= \frac{\bar{\gamma}U^*}{\alpha^* S^*} \omega_2, \quad \omega_{13}= \frac{\bar{\gamma}U^*}{\alpha^* S^*}\omega_3,\quad \omega_{14}=0,\;\;\textrm{with} \;\;\omega_2>0,\quad \omega_3>0.
\end{aligned}
\end{equation}
In a similar way, the left eigenvector of $J|_{E_0}$ can be denoted as  $\textbf{v}=(v_1,v_2,\dots,v_{14})$:
\begin{equation}\label{E:lvec}
\begin{aligned}
    v_1&=0,\;\; v_4=0,\;\; v_5=0,\;\; v_6=0,\;\; v_7=0,\;\;v_8=\frac{\bar{\gamma}}{\bar{\delta}}v_3,\;\;v_9=\frac{\bar{\gamma}}{\bar{\delta}}v_2, \;
    v_{10}=0,\;v_{11}=0,\\ v_{12}&=\frac{\alpha^* S^*}{\kappa U^*}v_2,\;\;v_{13}=\frac{\alpha^* S^*}{\kappa U^*}v_3,\;\;v_{14}=\frac{\alpha^* S^*}{\kappa U^*}(v_2+v_3),\;\;\textrm{with} \;\;v_2>0,\quad v_3>0.
\end{aligned}
\end{equation}
where $\bar{\gamma}=\gamma+\mu,\; \bar{\delta}=\delta+\gamma+\mu$. Thus one obtains
\begin{align*}
    \text{\textbf{w}} \cdot \text{\textbf{v}} &= \left(1+\frac{\bar{\gamma}}{\kappa}\right) (v_2\omega_2+v_3 \omega_3)>0.
\end{align*}
Consider $a$ to be the bifurcation coefficient presented in \cite{van2002reproduction}. Rewriting the model \eqref{Eq:S}-\eqref{Eq:V12} in the form of $\dot{x_i}=F_i(x_i)$, $i=\{1,2,\dots,14\}$, this coefficient is given by
\begin{equation}
    a=\sum_{k,i,j=1}^{14} v_k \omega_i \omega_j \frac{\partial^2 F_k}{\partial x_i \partial x_j} (0,0).\label{E:abif}
\end{equation}
Using Eqs. \eqref{E:lvec} and \eqref{E:abif}, it becomes
\begin{equation}
    \begin{aligned}
        a&= v_2 \sum_{i,j=1}^{14} \omega_i \omega_j \frac{\partial^2 F_2}{\partial x_i \partial x_j}+ v_3  \sum_{i,j=1}^{14} \omega_i \omega_j \frac{\partial^2 F_3}{\partial x_i \partial x_j}+ \frac{\bar{\gamma}}{\bar{\delta}}v_3  \sum_{i,j=1}^{14} \omega_i \omega_j \frac{\partial^2 F_8}{\partial x_i \partial x_j}\\
        &\qquad + \frac{\bar{\gamma}}{\bar{\delta}}v_2 \sum_{i,j=1}^{14} \omega_i \omega_j \frac{\partial^2 F_9}{\partial x_i \partial x_j}+ \frac{\alpha^* S^*}{\kappa U^*}v_2  \sum_{i,j=1}^{14} \omega_i \omega_j \frac{\partial^2 F_{12}}{\partial x_i \partial x_j}+ \frac{\alpha^* S^*}{\kappa U^*}v_3  \sum_{i,j=1}^{14} \omega_i \omega_j \frac{\partial^2 F_{13}}{\partial x_i \partial x_j}\\
        &\qquad +\frac{\alpha^* S^*}{\kappa U^*}(v_2+v_3)\sum_{i,j=1}^{14} \omega_i \omega_j \frac{\partial^2 F_{14}}{\partial x_i \partial x_j},
    \end{aligned}
\end{equation}
where
\begin{align*}
    F_2&=\frac{\alpha^* x_1(x_{12}+x_{14})}{x_{11}+x_{12}+x_{13}+x_{14}} -\bar{\gamma} x_2,  \\
    F_3&=\frac{\alpha^* x_1 (x_{13}+x_{14})}{x_{11}+x_{12}+x_{13}+x_{14}} - \bar{\gamma} x_3,\\
    F_8&= \frac{\sigma x_6 (x_{13}+x_{14})}{x_{11}+x_{12}+x_{13}+x_{14}} - (\bar{\gamma}+\delta) x_8,\\
    F_9&= \frac{\sigma x_7 (x_{12}+x_{14})}{x_{11}+x_{12}+x_{13}+x_{14}} - (\bar{\gamma}+\delta) x_9,\\
    F_{12} &=  \frac{\beta \left[x_{11} (x_{2}+x_{9})-   x_{12} (x_{3}+x_{8})\right]}{x_1+x_2+x_3+x_4+x_5+x_6+x_7+x_8+x_9+x_{10}} -\kappa x_{12},\\
    F_{13}&= \frac{\beta \left[x_{11} (x_{3}+x_{8})-   x_{13} (x_{2}+x_{9})\right]}{x_1+x_2+x_3+x_4+x_5+x_6+x_7+x_8+x_9+x_{10}} -\kappa x_{13},\\
     F_{14}&=\frac{\beta \left[x_{12} (x_{3}+x_{8})+ x_{13} (x_{2}+x_{9})\right]}{x_1+x_2+x_3+x_4+x_5+x_6+x_7+x_8+x_9+x_{10}} -\kappa x_{14}.\\
\end{align*}
Using Eq. \eqref{E:rvec}, we eventually obtain the bifurcation coefficient evaluated at the disease free equilibrium as
\begin{equation}\label{E:aE0}
   a|_{E_0} = - A_1 (v_2\omega_2+v_3\omega_3)(\omega_2+\omega_3) + \frac{A_2}{\alpha^*}(v_2+v_3)\omega_2\omega_3,
\end{equation}
with $A_1=\frac{2 \bar{\gamma}}{S^*}\left(\frac{\beta}{\kappa}+\frac{\bar{\gamma}}{\mu}\right)$, $A_2=\frac{2\sigma\gamma\nu \bar{\gamma}^2}{\bar{\nu}\bar{\delta}\mu S^*}$.
Here we want to determine the sign of $a$. Since we have assume that $v_2>0, v_3>0, \omega_2>0 $ and $\omega_3>0$, the bifurcation coefficient $a$ is found positive if
\begin{equation}\label{Eq:alpc}
    \alpha^*< \frac{A_2 (v_2+v_3)\omega_2\omega_3 }{A_1 (v_2\omega_2+v_3\omega_3)(\omega_2+\omega_3)}=\alpha_c.
\end{equation}
 In the notation given in \cite{van2002reproduction} one also obtains the coefficient
\begin{equation}
    b=\sum_{k,i=1}^{14} \omega_k \omega_i \frac{\partial^2 f_k}{\partial x_i \partial \alpha}(0,0),
\end{equation}
evaluated at the disease free equilibrium leads to
\begin{equation}
    b|_{E_0} =\frac{\bar{\gamma}^2}{\alpha^* \mu} (\omega_2+\omega_3)^2+\frac{\bar{\gamma}}{\alpha^*} (\omega_2^2+\omega_3^2)>0.
\end{equation}
In this context Eq. \eqref{E:aE0} can be rewritten as
\begin{equation}
     a|_{E_0} =-A_1(v_2\omega_2^2+v_3\omega_3^2)+(\frac{A_2}{\alpha^*}-A_1)(v_2+v_3)\omega_2\omega_3. \label{E:AE02}
\end{equation}
 The model exhibits backward bifurcation whenever $\alpha^*<\alpha_c$ in Eq. \eqref{Eq:alpc}.  Furthermore, it should be noted from Eq. \eqref{E:AE02} that the conditions ($a<0$) for a backward bifurcation in \cite{van2002reproduction} can not be satisfied if $\alpha>\frac{A_2}{A_1}$.
\begin{table}[H]
\centering
\caption{The values of the parameters (monthly units) in our model \ref{Eq:Model} are based on information about dengue epidemiology.}
\begin{tabular}{lllc}
\hline
\textbf{Parameter} & \textbf{Value (month$^{-1}$)} & \textbf{Notes} & \textbf{Reference} \\
\hline
$\Lambda_N$ & $12.8$ & Recruitment rate of humans & \cite{ferguson1999} \\ 
$\Lambda_M$ & $10^5$ & Recruitment rate of vectors, and $28,000$ per week& \cite{khan2014estimating} \\
$\mu$ & $0.00128$ & Avg. human lifespan $\approx 65-70$ yrs & \cite{ferguson1999transmission} \\
$\gamma$ & $2.0$& Avg. infectious period $\approx 15$ days, and $\gamma^{-1}=0.5$ month & \cite{ferguson1999transmission,esteva1998analysis} \\
$\nu$ & $0.111$ & Cross-immunity period & \cite{ferguson1999transmission, aguiar2008epidemiology, paz2021dengue,reich2013interactions} \\
$\delta$ & $0.01$ & Disease-induced mortality (rare, severe DHF/DSS) & \cite{halstead2003neutralization,sabchareon2012protective} \\ 
$\kappa$ & $1$ & Vector lifespan, and $\kappa \in [0.67, 1.0]$ & \cite{esteva1998analysis,sophia2025dengue} \\
$\alpha$ & $0.39$ & Primary infection rate from infected vectors  & \cite{ogunlade2023quantifying,esteva1998analysis,nishiura2006mathematical} \\ 
$\sigma$ & $0.1$–$0.8$ & Secondary infection rate (ADE scenarios) & \cite{ferguson1999transmission,esteva1998analysis,nishiura2006mathematical} \\
$\beta$ & $6$ & Vector infection rate from humans & \cite{esteva1998analysis, thongsripong2020investigation} \\
\hline
\end{tabular}
\label{tab:params_values}
\end{table}
Baseline parameter values in per-month units, assuming $\delta=0$ (constant human population) and quasi-constant vector population ($\Lambda_N=\mu N$, $\Lambda_M=\kappa M$).
\section{Numerical investigation of the two-strain dengue epidemics model}
\label{Sec:NumSim}

In this section, we perform a numerical investigation of the dynamical behavior of the two-strain dengue model \eqref{Eq:Model}, with emphasis on how qualitative outcomes change under variations of key parameters. The computations combine direct time integration of the ODE system with path-following (continuation) methods implemented in the software package COCO \cite{dankowicz2013}. COCO is a MATLAB-based platform for numerical continuation and bifurcation analysis of parametrized smooth ODEs; in spirit and functionality, it is comparable to established tools such as AUTO \cite{auto97} and MATCONT \cite{matcont}. Throughout, we keep the model structure and parameter meanings as defined previously (Sections~\ref{Sec:Model} and \ref{Sec:Prelim}), and use the baseline parameter values summarized in Table~\ref{tab:params_values}, unless explicitly stated otherwise.

\begin{figure}[H]
	\centering\psfrag{t}{\large$t$ \scriptsize[months]}\psfrag{S}{$S(t)$}\psfrag{Is}{$I_{1}(t)$, $I_{2}(t)$}\psfrag{Isec}{$(I_{12}+I_{21})(t)$}
	\includegraphics[width=\textwidth]{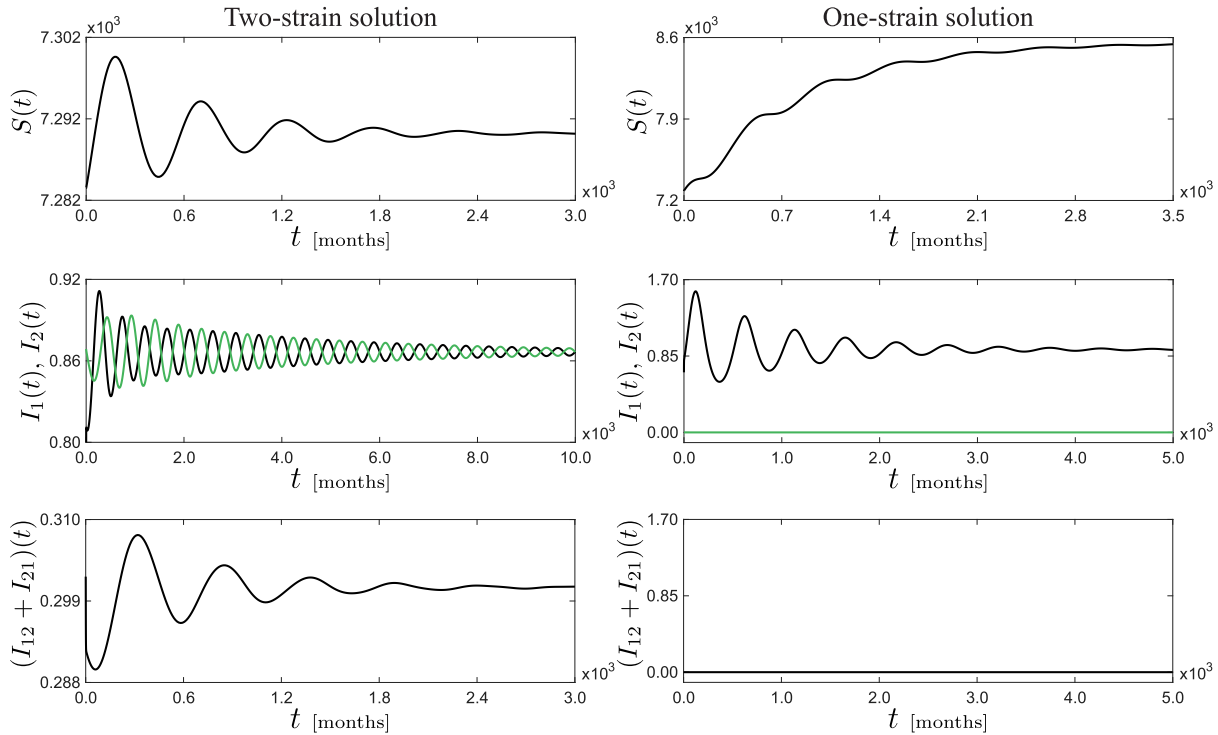}
	\caption{One- and two-strain transient responses of the dengue epidemics model \eqref{Eq:Model}, computed for the parameter values given in Table \ref{tab:params_values}, with $\sigma=0.45$. In this picture, the time histories for $I_{1}(t)$ and $I_{2}(t)$ are plotted in black and green, respectively.}\label{fig-sol-ini}
\end{figure}

To initiate our analysis, we adopt the parameter values of the two-strain dengue epidemics model \eqref{Eq:Model} as listed in Table \ref{tab:params_values}, with $\sigma=0.45$. With these values, the system response, obtained via direct numerical integration, is presented in Fig.\ \ref{fig-sol-ini}. The figure shows time series for the following state variables: $S(t)$ (susceptible individuals that can be infected by all dengue strains), $I_{1}$ (individuals infected with strain 1), $I_{2}$ (individuals infected with strain 2), and $I_{12}+I_{21}$ (individuals suffering secondary infections, i.e.\ those that have been infected by both strains, one after another). Interestingly, for the same set of parameter values, two types of behavior can be detected: one showing existence of both strains in the population, and one response showing one strain only ($I_{1}$ in this case). The first solution is obtained using nonzero initial conditions for all compartments, while the second scenario is produced by setting $I_{2}=R_{2}=S_{2}=I_{12}=I_{21}=R=V_{2}=V_{12}=0$.

Due to symmetry, another solution may be computed, corresponding to the existence of strain 2 only, but for the sake of brevity only the former case will be considered. The existence of only one strain can be clearly verified from the second column of Fig.\ \ref{fig-sol-ini}, where the solutions corresponding to $I_{2}$ (in green) and secondary infections $I_{12}+I_{21}$ are zero. As can be seen in the time histories, the system in both cases eventually settles down to an equilibrium, one with both strains present, and one with strain 1 only. These equilibria will serve as the baseline for our numerical exploration, where we will focus on how changes in critical parameters affect the steady state, with special attention put on variations of dengue infections with respect to parameter perturbations.

\begin{figure}[H]
	\centering
	\psfrag{a}{\large(a)}\psfrag{b}{\large(b)}\psfrag{t}{\large$t$ \scriptsize[months]}\psfrag{S}{$S(t)$}\psfrag{Is}{$I_{1}(t)$, $I_{2}(t)$}\psfrag{Vs}{$V_{1}(t)$, $V_{2}(t)$}\psfrag{al}{\large$\alpha$}\psfrag{sirat}{\large$\sigma/\alpha$}\psfrag{nu}{\large$\nu$}\psfrag{It}{\large$I_{\text{\tiny Tot}}$}\psfrag{Isec}{\large$I_{\text{\tiny Sec}}$}
	\includegraphics[width=0.95\textwidth]{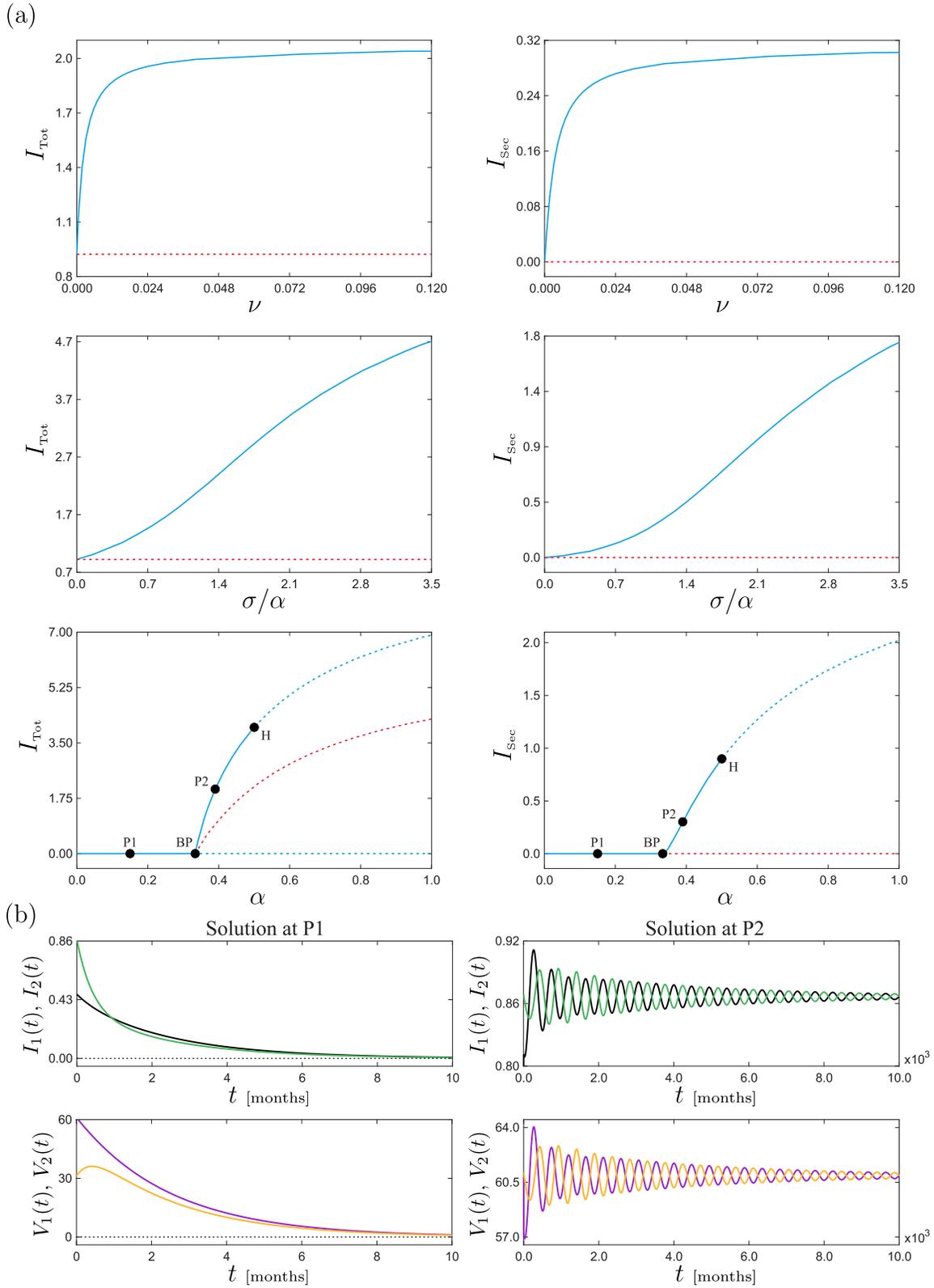}
	\caption{(a) One-parameter continuation of the two-strain (blue curve) and one-strain (red curve) equilibria calculated from Fig.\ \ref{fig-sol-ini}. In these diagrams, the vertical axes are given by $I_{\text{\tiny Tot}}=I_{1}+I_{2}+I_{12}+I_{21}$ and $I_{\text{\tiny Sec}}=I_{12}+I_{21}$. The BP and H points stand for branching ($\alpha\approx0.33355$) and Hopf ($\alpha\approx0.50012$) bifurcations, respectively. Solution branches of stable equilibria are plotted with solid lines, while dashed lines denote instability. (b) System responses calculated at the test points P1 ($\alpha=0.15$) and P2 ($\alpha=0.39$). Here, the color code is as follows: $I_{1}(t)$ (black), $I_{2}(t)$ (green), $V_{1}(t)$ (purple) and $V_{2}(t)$ (yellow).}\label{fig-bif-diags}
\end{figure}

We now apply numerical continuation (path-following) methods to investigate how the equilibria identified earlier evolve under gradual variation of selected parameters. The outcome of this procedure is illustrated in Fig.\ \ref{fig-bif-diags}(a). Here, two solution measures are displayed: $I_{\text{\tiny Tot}}=I_{1}+I_{2}+I_{12}+I_{21}$, showing total number of dengue infections (primary and secondary), and $I_{\text{\tiny Sec}}=I_{12}+I_{21}$, measuring secondary infections only. Continuation of two-strain and one-strain equilibria are plotted using blue and red colors, respectively. The first parameter considered for the numerical continuation is the cross-immunity factor $\nu$, with $1/\nu$ being the average cross-immunity duration. As explained earlier, this parameter defines how long an individual remains immune to a secondary dengue infection. In the absence of strain 2 (see second column in Fig.\ \ref{fig-sol-ini}), this parameter has no effect in the system, due to which a horizontal straight line (red) is obtained in Fig.\ \ref{fig-bif-diags}(a). When both strains are present, on the contrary, $\nu$ does affect the level of dengue infections, as can be seen from the behavior of the blue line. As can be expected, larger values of $\nu$ (meaning shorter cross-immunity duration) increase the level of both total and secondary dengue infections.

Another parameter considered in our study is the secondary infection rate $\sigma$. For a better understanding, in Fig.\ \ref{fig-bif-diags}(a) we have displayed the ratio $\sigma/\alpha$ (with $\alpha$ fixed) to compare the levels of secondary and primary infection rates. Analogously to the case studied before, this ratio has no effect whatsoever when only one strain is present in the system. However, the situation clearly changes when both strains are circulating in the model, as can be seen in Fig.\ \ref{fig-bif-diags}(a). As before, when $\sigma$ increases higher levels of total and secondary dengue infections are produced in the model. Next, we will consider variations in the primary infection rate $\alpha$. As can be seen in Fig.\ \ref{fig-bif-diags}(a), for low values of this parameter a stable disease-free equilibrium is present, with basic reproduction number $R_{0}<1$. As $\alpha$ grows, $R_{0}$ grows too (see Section \ref{Subsec:R0}), until the critical value 1 is reached, for $\alpha\approx0.33355$. Here, a classical forward bifurcation occurs, giving rise to three branches: unstable disease-free, stable two-strain and unstable one-strain equilibria, see Fig.\ \ref{fig-bif-diags}(a). Notice that the red curve corresponding to one-strain solutions shows increasing behavior of the total number of dengue infections (influenced by the primary infection $I_{1}$), while secondary infections remain zero, as expected. On the other hand, the two-strain solution shows growing behavior in both dengue indicators. To confirm the numerical observations, two sample solutions are calculated via direct numerical integration in Fig.\ \ref{fig-bif-diags}(b), corresponding to the test points P1 (disease-free) and P2 (two-strain solution).

\begin{figure}[H]
	\centering
	\psfrag{a}{\large(a)}\psfrag{b}{\large(b)}\psfrag{t}{\large$t$ \scriptsize[months]}\psfrag{Isecs}{$I_{12}(t)$, $I_{21}(t)$}\psfrag{S}{$S$}\psfrag{U}{$U$}\psfrag{Is}{$I_{1}(t)$, $I_{2}(t)$}\psfrag{Ss}{$S_{1}(t)$, $S_{2}(t)$}\psfrag{Vs}{$V_{1}(t)$, $V_{2}(t)$}
	\includegraphics[width=0.9\textwidth]{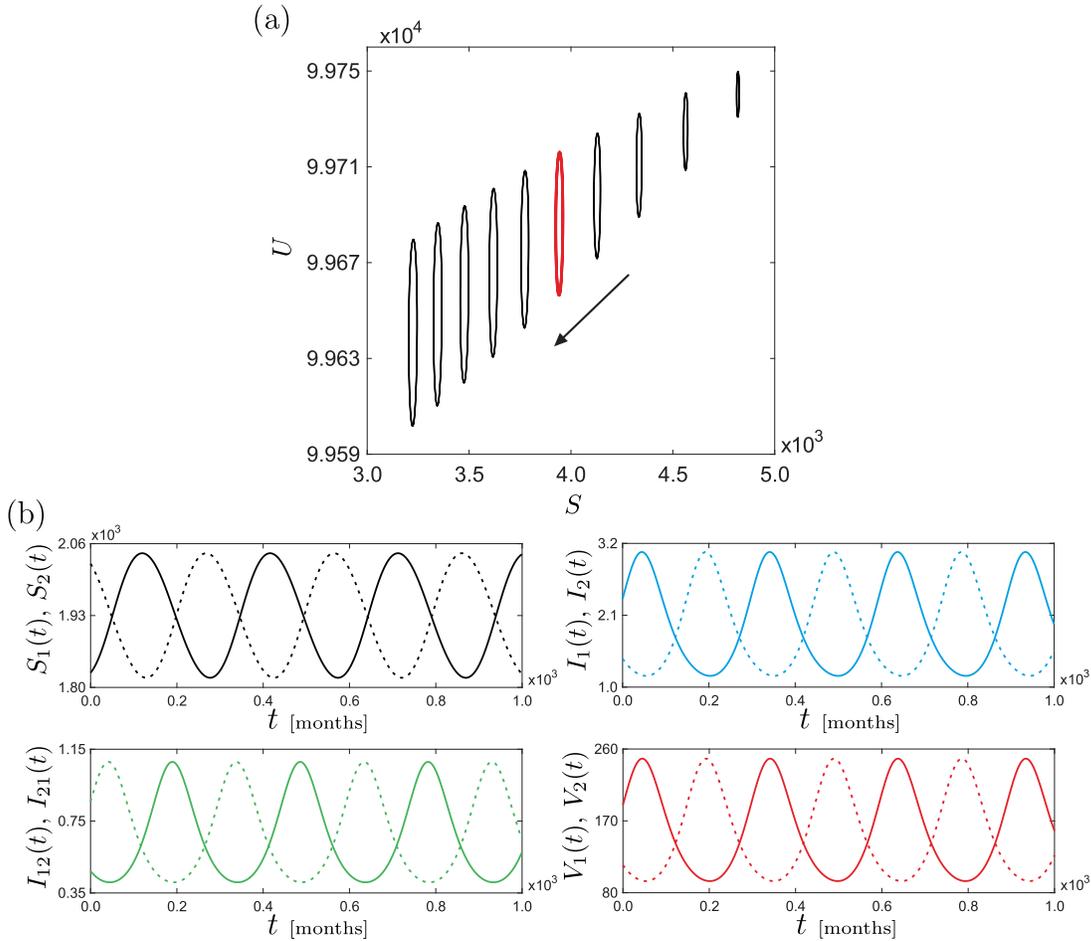}
	\caption{(a) Family of periodic solutions produced by the Hopf bifurcation detected in Fig.\ \ref{fig-bif-diags}(a), calculated from $\alpha=0.529$ to $\alpha=0.75$. The arrow indicates the direction of increasing $\alpha$. (b) Time histories corresponding to the solution highlighted in red ($\alpha=0.627$) in panel (a). In these plots, solid and dashed lines correspond to state variables with leading subscript 1 and 2, respectively.}\label{fig-per-orbs}
\end{figure}

Our numerical investigation also reveals the presence of a Hopf bifurcation, detected when the primary infection rate reaches the critical value $\alpha\approx0.50012$. Here, the two-strain equilibrium becomes unstable (represented by a dashed line), and a family of stable periodic
solutions is born, see Fig.\ \ref{fig-per-orbs}. In this figure, panel (a) displays a sequence of periodic orbits with growing amplitude as $\alpha$ increases beyond the Hopf point (from $\alpha=0.529$ to $\alpha=0.75$). A test solution is depicted in panel (b), showing the oscillatory behavior of both human and vector compartments.

\section{Conclusion}
\label{Sec:Conc}

In this work, we proposed  a novel, deterministic, two-strain host-vector dengue transmission model.
The model advances previous two-strain frameworks by distinguishing between primary and secondary human infections, integrating ADE effects through distinct transmission parameters $(\alpha,\sigma)$, and allowing vectors to host both single and double infections, thereby providing a more realistic representation of dengue ecology.

Analytically, we derived the basic reproduction number using the next-generation matrix method, obtaining
$ R_0 = \sqrt{\frac{\alpha \beta}{\kappa (\gamma + \mu)}},$
and established that the disease-free equilibrium is locally asymptotically stable whenever $R_0 < 1$ and unstable for $R_0 > 1$. The one-strain endemic equilibrium loses stability when $\sigma > \alpha$, illustrating how ADE enables the invasion of a heterologous strain. Through center-manifold theory, we demonstrated that the system admits backward bifurcation under suitable parameter regimes, implying that dengue persistence is possible even when $R_0 < 1$, due to the coexistence of stable disease-free and endemic equilibria. This mathematical property highlights the challenge of dengue eradication: simple threshold control strategies based solely on $R_0$ may fail when enhancement mechanisms are strong.

Our numerical continuation and bifurcation analysis (performed with COCO) revealed a rich dynamical landscape consistent with the analytical findings. For baseline parameters, the model exhibits bistability between single- and two-strain endemic equilibria. As the primary infection rate $\alpha$ increases, a forward (branching) bifurcation occurs near $\alpha \approx 0.33355$, marking the transition from the disease-free to the endemic regime. A subsequent Hopf bifurcation at $\alpha \approx 0.50012$ gives rise to a family of stable periodic orbits, characterized by oscillations in both human and vector infection classes. The amplitude of these oscillations grows with $\alpha$, indicating how ADE and high transmission intensity can produce recurrent epidemic outbreaks. Variations in cross-immunity loss $\nu$ and the secondary-to-primary infection ratio $\sigma/\alpha$ further demonstrated how shorter immunity durations and stronger enhancement intensify total and secondary infection levels, respectively. Together, these results confirm that ADE, temporary immunity, and vector co-infection jointly generate multistability and oscillatory behavior observed in multi-strain dengue epidemics.

From an epidemiological perspective, these findings suggest that reducing $R_0$ below unity may be insufficient for disease elimination when backward bifurcation or bistability occur. Control measures that reduce secondary transmission ($\sigma$), prolong temporary immunity ($1/\nu$), or mitigate ADE effects could play a more effective role than focusing solely on primary infection reduction. Moreover, understanding how parameter interactions drive oscillatory outbreaks can inform vaccination design, especially for vaccines that modulate immune enhancement or cross-protection.

Future work could extend the present framework by incorporating age or spatial structure, stochastic effects, or seasonally forced mosquito dynamics to better capture the temporal variability observed in real-world dengue outbreaks. The model can also be coupled with optimal control formulations to assess the efficiency of combined interventions such as vaccination, vector control, and public health campaigns. Ultimately, the proposed two-strain host–vector model provides a rigorous and biologically grounded platform for exploring complex dengue transmission dynamics and for guiding the development of evidence-based control strategies in dengue-endemic regions.

\newpage
\section{Appendix: Linearization}\label{Sec:Appendix1}
Linearization around a generic steady state
\[E= \left(S^{**}, I_1^{**}, I_2^{**}, R_1^{**}, R_2^{**}, S_1^{**}, S_2^{**}, I_{12}^{**}, I_{21}^{**}, R^{**},U^{**}, V_1^{**}, V_2^{**}, V_{12}^{**}\right)\]
 leads to the Jacobian matrix
\begin{equation}
    \mathcal{J}=\left[J_{ij}\right]_{14 \times 14}, \qquad \textrm{for}\; i,j=1,2,\dots,14,
\end{equation}
where Jacobian can be re-expressed as
$$\mathcal{J}=[\mathcal{M}_{{14}\times7} \;\;\mathcal{N}_{{14}\times7}]$$
with

\[
\mathcal{M}=
\resizebox{\textwidth}{!}{
$\begin{bmatrix}
-{\alpha M^{**} A_{V}} -\mu & 0 & 0 & 0 & 0 & 0 & 0 \\

 {\alpha M^{**} A_1} & -\bar{\gamma} & 0 & 0 & 0 & 0 & 0 \\

 {\alpha}{M^{**}} A_2 & 0 & -\bar{\gamma} & 0 & 0 & 0 & 0 \\

 0 & \gamma & 0 & -\bar{\nu} & 0 & 0 & 0 \\

 0 & 0 & \gamma & 0 & -\bar{\nu} & 0 & 0 \\

 0 & 0 & 0 & \nu & 0 & -{\sigma}{M^{**}} A_2 -\mu & 0 \\

 0 & 0 & 0 & 0 & \nu & 0 & -{\sigma}{M^{**}}A_1-\mu \\

 0 & 0 & 0 & 0 & 0 & {\sigma}{M^{**}} A_2 & 0 \\

 0 & 0 & 0 & 0 & 0 & 0 & {\sigma}{M^{**}} A_1 \\

 0 & 0 & 0 & 0 & 0 & 0 & 0 \\

  \beta U^{**} A_I & -\beta U^{**} A_N& -\beta U^{**} A_N & \beta U^{**} A_I &\beta U^{**} A_I & \beta U^{**} A_I & \beta U^{**} A_I  \\

A_{UV_1} &  \beta U^{**}/N^{**} +A_{UV_1} & -\beta V^{**}_1/N^{**}+A_{UV_1} & A_{UV_1} & A_{UV_1} & A_{UV_1} & A_{UV_1} \\

A_{UV_2} & -\beta V^{**}_1/N^{**}+A_{UV_2} & \beta U^{**}/N^{**}+A_{UV_2}  & A_{UV_2}& A_{UV_2} &A_{UV_2} & A_{UV_2}  \\

A_{VV} & \beta V^{**}_2/N^{**} +A_{VV}  & \beta V^{**}_1/N^{**} +A_{VV}   & A_{VV} & A_{VV} & A_{VV} & A_{VV}\\
\end{bmatrix}$.
}
\]

and
\[
\mathcal{N}=
\resizebox{\textwidth}{!}{
$\begin{bmatrix}
0 & 0 & 0 & {\alpha S^{**} A_{V}} & -\alpha S^{**} A_U & -\alpha S^{**} A_U & -\alpha S^{**} A_{UV}\\

0 & 0 & 0 & - {\alpha S^{**} A_1} & {\alpha S^{**} A_{U2}} & - {\alpha S^{**} A_1} & {\alpha S^{**} A_{U2}}\\

 0 & 0 & 0 & -\alpha S^{**} A_{2} & -\alpha S^{**} A_{2} & \alpha S^{**} A_{U1} & \alpha S^{**} A_{U1}\\

 0 & 0 & 0 & 0 & 0 & 0 & 0\\

 0 & 0 & 0 & 0 & 0 & 0 & 0\\

0 & 0 & 0 & \sigma S^{**}_1 A_2 &  \sigma S^{**}_1 A_2 & -\sigma S^{**}_1 A_{U1} & -\sigma S^{**}_1 A_{U1}\\

 0 & 0 & 0 & \sigma S^{**}_2 A_1  & -\sigma S^{**}_2 A_{U2} & \sigma S^{**}_2 A_1 & -\sigma S^{**}_2 A_{U2}\\

 -\bar{\delta} & 0 & 0 & -\sigma S^{**}_1 A_2 & -\sigma S^{**}_1 A_2 & {\sigma}{S^{**}_1} A_{U1}  & {\sigma}{S^{**}_1} A_{U1} \\

 0 & -\bar{\delta} & 0 & -{\sigma}{S^{**}_2} A_1  & {\sigma}{S^{**}_2} A_{U2} & -{\sigma}{S^{**}_2} A_1  & {\sigma}{S^{**}_2} A_{U2}\\

  \gamma & \gamma & -\mu & 0 & 0 & 0 & 0\\

  -\beta U^{**} A_N  &  -\beta U^{**} A_N  &  \beta U^{**} A_I & -{\beta}{N^{**}}A_I-\kappa & 0 & 0 & 0\\

 -\beta V^{**}_1/N^{**}+A_{UV_1} & \beta U^{**}/N^{**} +A_{UV_1} & A_{UV_1} & N^{**} A_{I1} & -N^{**} A_{I2}-\kappa  & 0 & 0\\

  \beta U^{**}/N^{**}+A_{UV_2}& -\beta V^{**}_2/N^{**}+A_{UV_2}  & 0 & N^{**} A_{I2} & 0 & -N^{**} A_{I1} -\kappa & 0 \\

  \beta V^{**}_1/N^{**}+A_{UV_2}   & \beta V^{**}_2/N^{**} +A_{UV_2}   & A_{UV_2} & 0 & N^{**} A_{I2} & N^{**} A_{I1} & -\kappa\\
\end{bmatrix}$.
}
\]

where $\bar{\gamma}=\gamma+\mu$, $\bar{\nu}=\nu+\mu$, $\bar{\delta}=\gamma+\mu+\delta$ and

\begin{align}
    A_{V}&=\frac{V^{**}_1+V^{**}_2+2 V^{**}_{12}}{(M^{**})^2},\;
    A_{1}=\frac{V^{**}_1+V^{**}_{12}}{(M^{**})^2},\;
    A_{2}=\frac{V^{**}_2+V^{**}_{12}}{(M^{**})^2},\;
    A_{U}=  \frac{(U^{**}-V^{**}_{12})}{(M^{**})^2},\; \nonumber \\
    A_{UV} &= \frac{V^{**}_1 + V^{**}_2+2U}{(M^{**})^2},
    A_{U1} = \frac{V^{**}_1+U^{**}}{(M^{**})^2},\; A_{U2} = \frac{V^{**}_2+U^{**}}{(M^{**})^2},
    A_I = \frac{I^{**}_1+I^{**}_2+I^{**}_{12} + I^{**}_{21}}{(N^{**})^2}, \nonumber\\
    A_N&=\frac{S^{**}+R^{**}_1+R^{**}_2+S^{**}_1+S^{**}_2+R^{**}}{(N^{**})^2},
    A_{I1} =\beta \frac{I^{**}_1+I^{**}_{21}}{(N^{**})^2},\;A_{I2} =\beta \frac{I^{**}_2+I^{**}_{12}}{(N^{**})^2}, \; \nonumber\\ A_{UV_1} &= -U^{**} A_{I1}+V^{**}_1 A_{I2},\; A_{UV_2}= - U^{**} A_{I2}+V^{**}_2 A_{I1},
    A_{VV}  = -V^{**}_1 A_{I2} -V^{**}_2 A_{I1}.
\end{align}

\color{black}

\bibliographystyle{ieeetr}
\bibliography{twostrainbg.bst}

\end{document}